\documentclass[12pt]{amsart}
\usepackage{amscd}

\makeatletter
\def\pf{%
  \par\topsep6\p@\@plus6\p@
  \trivlist
  \item[\hskip\labelsep\it\proofname.]\ignorespaces}

\makeatother

\let\:=\colon
\def\i{^{-1}}
\def\CC{{\mathbf C}}
\def\QQ{{\mathbf Q}}
\def\PP{{\mathbf P}}

\def\F{{\mathcal F}}
\def\G{{\mathcal G}}

\def\I{{\mathcal I}}
\def\E{{\mathcal E}}

\def\Z{{\mathcal Z}}
\def\OO{{\mathcal O}}

\def\L{{\mathcal L}}
\def\x{\times}
\def\*{\otimes}
\def\v{^{\vee}}
\def\iso{\simeq}

\def\sub{\subseteq}
\def\Ext{\operatorname{Ext}}
\def\Hilb{\operatorname{Hilb}}
\def\Pic{\operatorname{Pic}}
\def\dsum{\oplus}

\newtheorem{thm}{Theorem}[section]

\newtheorem{lem}[thm]{Lemma}
\newtheorem{prop}[thm]{Proposition}

\theoremstyle{definition}
\newtheorem{defn}[thm]{Definition}

\theoremstyle{remark}

\newtheorem{notation}{Notation}

\newtheorem{ack}{Acknowledgements}

\numberwithin{equation}{section}

\newcommand{\thmref}[1]{theorem~\ref{#1}}

\newcommand{\lemref}[1]{lemma~\ref{#1}}

\begin{document}
\bibliographystyle{plain}

\title
{ Some Donaldson invariants of $\CC\PP^2$}

\author{Geir Ellingsrud}
\address{Mathematical Institute\\
         University of Oslo\\P.~O.~Box~1053\\
         N--0316 Oslo, Norway}
         \email{ellingsr@@math.uio.no}
\author{Joseph Le Potier}
\address{Universit\'e Paris 7, UFR de
         Math\'ematiques et Institut de Math\'e\-ma\-tiques de
         Jussieu,
         Case Postale 7012, 2, Place Jussieu,
         F--75251 Paris Cedex 05}
         \email{jlp@@mathp7.jussieu.fr}
\author{Stein A.~Str{\o}mme}
\address{Mathematical Institute\\
         University of Bergen\\All\'eg 55\\
         N--5007 Bergen, Norway}
         \email{stromme@@mi.uib.no}

\date{\today}


\maketitle

\begin{center}
\emph{In memory of the victims of the Kobe earthquake}
\end{center}

\section*{Introduction}
For an integer $n\ge2$, let $q_{4n-3}$ be the coefficient of
the Donaldson polynomial of degree $4n-3$ of
$P=\CC\PP^2$.  An interpretation of $q_{4n-3}$ in an
algebro-geometric context is the following.  Let $M_n$
denote the Gieseker-Maruyama moduli space of semistable
coherent sheaves on $P$ with rank 2 and Chern classes
$c_1=0$ and $c_2=n$.  For such a sheaf $F$, the
Grauert-M\"ulich theorem implies that the
restriction of $F$ to a general line $L\sub P$ splits as
$F_L \iso \OO_L\dsum\OO_L$, and that the exceptional lines
form a curve $J(F)$ of degree $n$ in the dual
projective plane $P\v$.  The association $F\mapsto J(F)$
is induced from a morphism of algebraic varieties,
called the Barth map, $f_n\: M_n \to P_n$.  Here
$P_n=\PP^{n(n+3)/2}$ is the linear system parameterizing all
curves of degree $n$ in $P\v$. Let $H\in\Pic(P_n)$ be the
hyperplane class and let $\alpha = f_n^*H$.
The interpretation of
the Donaldson invariant is:
\[
q_{4n-3} = \int_{M_n} \alpha^{4n-3}.
\]
Thus $q_{4n-3}$ is the degree of $f_n$ times the degree of
its image.  From \cite{Bart-2} it follows that $f_n$ is
generically finite for all $n\ge2$, that $f_2$ is an
isomorphism and $q_5=1$, and that $f_3$ is of degree 3 and
$q_9=3$.  Le Potier \cite{LePo} proved that $f_4$ is
birational onto its image
and that $q_{13}=54$.  The value of $q_{13}$ has
also been computed independently by Tikhomirov and Tyurin
\cite[prop.~4.1]{Tyur-1}
and by Li and Qin \cite[thm.~6.29]{Li-Qin}.

The main result in the present note is the following
\begin{thm} \label{thm1}
   $q_{17}=2540$ and $q_{21}=233208$.
\end{thm}

The proof consists of two parts.  The first part, treated in
this note, is to express $q_{4n-3}$ in terms of certain
classes on the Hilbert scheme of length-$(n+1)$ subschemes
of $P$.  This is theorems \ref{thm2} and \ref{thm3} below.
The second part is to evaluate these classes numerically.
This has been carried out in \cite[prop.~4.2]{Elli-Stro-5}.


Let $H_{n+1}=\Hilb^{n+1}_P$ denote the Hilbert scheme
parameterizing closed subschemes of $P$ of length $n+1$.
There is a universal closed subscheme $\Z\sub H_{n+1}\x P$.
Consider the vector bundles
\[
\E = R^1{p_1}_* (\I_{\Z}\*{p_2}^* \OO_{P}(-1))\text{ and }
\G = R^1{p_1}_* \I_{\Z}
\]
on $H_{n+1}$ of ranks $n+1$ and $n$, respectively, and
the linebundle
\[\L = \det(\G) \* \det(\E)\i.\]

\begin{thm} \label{thm2}
Let the notation be as above. Then
\[ q_{17} = \int_{H_6} s_{12}(\E\*\L) \quad\text{and}\quad
q_{21} = \dfrac25 \int_{H_7} s_{14}(\E\*\L).
\]
\end{thm}

This result was obtained both by Tikhomirov and Tyurin
\cite{Tyur-Tikh}, using the method of ``geometric
approximation procedure'' and by Le Potier \cite{LePo-3},
using ``coherent systems''. We present in this note what we
believe is a considerably simplified proof, which is
strongly hinted at on the last few pages of
\cite{Tyur-Tikh}.

The formula for $q_{17}$ is a special case of the
following formula:

\begin{thm}\label{thm3}
For $2\le n\le 5$, we have
\[q_{4n-3} =
\dfrac1{2^{5-n}}\int_{H_{n+1}} c_1(\L)^{5-n} s_{3n-3}(\E\*\L).
\]
\end{thm}
With this it is also easy to recompute $q_5$, $q_9$, and
$q_{13}$ using similar techniques as in \cite{Elli-Stro-5}.

\begin{notation}
We let $h$, $h\v$, and $H$ be the hyperplane classes in
$P$, $P\v$, and $P_n$, respectively. In general, if
$\omega$ is a divisor class, we denote by $\OO(\omega)$
the corresponding linebundle and its natural pullbacks.
\end{notation}

\begin{ack} This work is heavily inspired by conversations
with A.~Tyurin, and we thank him for generously sharing his
ideas.   We would also like to express our gratitude
towards the Taniguchi Foundation.
\end{ack}

\section{Hulsbergen sheaves}
Barth \cite{Bart-2} used the term Hulsbergen bundle to
denote a stable rank-2 vector bundle $F$ on $P$ with
$c_1(F)=0$ and $H^0(P,F(1))\ne0$.  We modify this definition
a little as follows:

\begin{defn} A \emph{Hulsbergen sheaf} is a coherent
sheaf $F$ on $P$ which admits a non-split short exact
sequence (\emph{Hulsbergen sequence})
\begin{equation} \label{Hulsbergen}
0 \to \OO_P \to F(1) \to \I_Z(2) \to 0,
\end{equation}
where $Z\sub P$ is a closed subscheme of finite length
(equal to $c_2(F)+1$).
\end{defn}

Note that a Hulsbergen sheaf is not necessarily semistable
or locally free.  However:

\begin{lem}\label{GM}
Let $F$ be a Hulsbergen sheaf with $c_2(F)=n>0$.  Then the
set $J(F)\sub P\v$ of exceptional lines for $F$ is a curve
of degree $n$, defined by the determinant of the bundle map
\[
m\: H^1(P,F(-2))\*\OO_{P\v}(-1) \to H^1(P,F(-1))\*\OO_{P\v}
\]
induced by multiplication with a variable linear form.
\end{lem}

\begin{pf}
First note from the Hulsbergen sequence that the two
co\-ho\-mo\-logy groups have dimension $n$.  It is easy to see
that any Hulsbergen sheaf is slope semistable, in the sense
that it does not contain any rank-1 subsheaf with positive
first Chern class.  Thus by \cite[thm.~1]{Bart-1}, $F_L \iso
\OO_L \dsum \OO_L$ for a general line $L$.  On the other
hand, it is clear that a line $L$ is exceptional if and only
if $m$ is not an isomorphism at the point $[L]\in P\v$.
\end{pf}

It is straightforward to construct a moduli space for
Hulsbergen sequences.  For any length-$(n+1)$ subscheme
$Z\sub P$, the isomorphism classes of extensions
\eqref{Hulsbergen} are parameterized by
$\PP(\Ext^1_P(\I_Z(2),\OO_P)\v)$.  By Serre duality,
\[
\Ext^1_P(\I_Z(2),\OO_P)\v \iso H^1(P,\I_Z(-1)).
\]
For varying $Z$, these vector spaces glue together to form
the vector bundle $\E$ over $H_{n+1}$, hence $D_n=\PP(\E)$
is the natural parameter space for Hulsbergen sequences.
Let $\OO(\tau)$ be the associated tautological quotient
linebundle.  For later use, note that for any divisor class
$\omega$ on $H_{n+1}$, we have
$\pi_*(\tau+\pi^*\omega)^{k+n} = s_k(\E(\omega))$, where
$\pi\: D_n \to H_{n+1}$ is the natural map \cite{IT}.

The tautological quotient $\pi^*\E \to \OO(\tau)$
gives rise to a short exact sequence on $D_n\x P$:
\[
0 \to \OO(\tau) \to \F(h) \to (\pi\x1)^*\I_{\Z}(2h) \to 0
\]
which
defines a complete family $\F$ of Hulsbergen sheaves.

As we noted earlier, a Hulsbergen sheaf is not necessarily
semistable.  On the other hand, the \emph{generic}
Hulsbergen sheaf is stable if $n\ge 2$.  It follows
that the family $\F$ induces a \emph{rational} map $g_n\:
D_n \to M_n$.  By \lemref{GM} above, there is also a Barth
map $b_n\: D_n \to P_n$, defined everywhere, and by
construction, the following diagram commutes:
\begin{equation}
\begin{CD}
             D_n  @>b_n>> P_n \\
             @V{g_n}VV    @VV{||}V \\
             M_n  @>>f_n> P_n
\end{CD}
\end{equation}
\begin{prop} Put $\lambda=c_1(\pi^*\L)$. Then
$b_n^*H = \tau+\lambda$.
\end{prop}

\begin{pf}  Let $L\sub P$ be a line. Twist the
universal Hulsbergen sequence by $-2h$ and $-3h$ respectively.
Multiplication by an equation for $L$ gives rise to the
vertical arrows in a commutative diagram with exact rows on
$D_n\x P$:
\[
\begin{CD}
0 @>>> \OO(\tau-3h) @>>> \F(-2h) @>>> (\pi\x1)^*\I_Z(-h)
@>>>0 \\
@. @VVV @VVV @VVV @.\\
0 @>>> \OO(\tau-2h) @>>> \F(-h) @>>> (\pi\x1)^*\I_Z @>>>0
\end{CD}
\]
Pushing this down via the first projection, we get the
following exact diagram on $D_n$:
\[\begin{CD}
0@>>>R^1{p_1}_* \F(-2h) @>>> \pi^*\E @>>> \OO(\tau)@>>> 0\\
@.    @Vm_LVV                  @VVV        @VVV\\
0@>>>     R^1{p_1}_* \F(-h) @>\iso>> \pi^*\G  @>>> 0
\end{CD}
\]
Here the last map of the top row is nothing but the
tautological quotient map on $\PP(\E)$.  Let $A(L)\sub D_n$
be the set of Hulsbergen sequences where $L$ is an
exceptional line for the middle term.  Clearly, $A(L)$ is
the degeneration locus of the left vertical map $m_L$ above.
Hence the divisor class of $A(L)$ is
\[
\begin{aligned}
[A(L)]&=  c_1(R^1{p_1}_* \F(-h)) - c_1(R^1{p_1}_* \F(-2h)) \\
      &=  \pi^*c_1(\G) - \pi^*c_1(\E) + \tau  \\
      &=  \tau+\lambda.
\end{aligned}
\]
On the other hand, $A(L)$ is the inverse image of a
hyperplane in $P_n$ under $b_n$, so its divisor class is
$b_n^*H$.
\end{pf}

\section{The case $n\le 5$}
\begin{prop}
For $2\le n\le 5$, the rational map $g_n$ is dominating, and
the general fiber is isomorphic to $\PP^{n-5}$.
For $n\ge 5$, the map $g_n$ is generically injective with
image of codimension $n-5$.  In particular, $g_5$ is
birational.
\end{prop}

\begin{pf}
Everything follows from the observation that the fiber over
a point $[F]\in M_n$ in the image of $g_n$ is the
projectivization of $H^0(P,F(1))$, and that for general such
$F$, this vector space has dimension
$h^0(F(1))=\max(1,6-n)$, which is easily seen from \eqref{Hulsbergen}.
  The assertion about the
codimension follows from a dimension count: $\dim(M_n)=4n-3$
and $\dim(D_n)=3n+2$.
\end{pf}

The first half of \thmref{thm2} now follows:
First of all, since $g_5$ is birational, the two morphisms
$f_5$ and $b_5$ have the same image and the same degree.
Therefore $q_{17}$ can be computed as
\[
q_{17} = \int_{D_5} H^{17} =\int_{D_5} (\tau+\lambda)^{17}
       = \int_{H_6} s_{12}(\E\*\L).
\]

For \thmref{thm3}, let $L_1,\dots,L_{5-n}$ be general lines
in $P$, and let $B_n\sub D_n$ be the locus of Hulsbergen
sequences where the closed subscheme $Z$ meets all these
$5-n$ lines.  The cohomology class of $B_n$ in
$H^*(D_n)$ is $\lambda^{5-n}$.

\begin{lem} \label{cover}
Let $2\le n\le5$.
The general nonempty fiber of $g_n$  meets $B_n$ in
$2^{5-n}$ points, hence the rational map
$g_n|_{B_n}\: B_n \to M_n$ is
 dominating and generically finite, of degree
$2^{5-n}$.
\end{lem}

\begin{pf}
The general nonempty fiber is of the form
$\PP(H^0(P,F(1))\v)$.  It suffices to show that the
restriction of $\L$ to this fiber has degree 2 (if $n<5$).
For this, it suffices to consider a linear pencil in the
fiber.  So let $\sigma_0$ and $\sigma_1$ be two independent
global sections of $F(1)$, and consider the pencil they
span.  Now $\sigma_0\wedge \sigma_1 \in
H^0(P,\wedge^2F)=H^0(P,\OO_P(2))$
is the equation of a conic $C\sub P$ which contains the zero
scheme $V(t_0\sigma_0 + t_1\sigma_1)$ of each section in the
pencil, $(t_0,t_1)\in\PP^1$.  Since $C$ meets a general line
in two points, it follows that there are exactly two members
of the pencil whose zero set meets a general line.
\end{pf}

To complete the proof of \thmref{thm3}, by \lemref{cover}
we now have for
$2\le n\le5$:
\[\begin{aligned}
2^{5-n}\,q_{4n-3} &= 2^{5-n}\int_{M_n} H^{4n-3} \\
        &=\int_{B_n} (\tau+\lambda)^{4n-3} \\
        &=\int_{D_n}
            \lambda^{5-n}\,(\tau+\lambda)^{4n-3} \\
        &=\int_{H_{n+1}}c_1(\L)^{5-n}
        \,s_{3n-3}(\E\*\L).
\end{aligned}
\]
This completes the proof of the theorems for $n\le 5$.

\section{The case $n=6$}
For $n\ge6$ the techniques above will say something about
the restriction of the Barth map to the Brill-Noether locus
$B\sub M_n$ of semistable sheaves whose first twist admit a
global section.  For general $n$ this locus is too small to
carry enough information about $M_n$, but in the special
case $n=6$, it is actually a divisor, whose divisor class
$\beta=[B]$ we can determine.  Now $\Pic(M_n)\*\QQ$ has rank
2, generated by $\alpha$ and $\delta=[\Delta]$, the class of
the locus $\Delta\sub M_n$ corresponding to non-locally free
sheaves
\cite{LePo-1}.

\begin{prop}
In $\Pic(M_6)\*\QQ$, the following relation holds:
\[
\beta = \frac52 \,\alpha - \frac12\,\delta.
\]
\end{prop}

\begin{pf}
Let $\xi\:X\to M_6$ be a morphism induced by a flat family
$\F$ of semistable sheaves on $P$, parameterized by some
variety $X$.  For certain divisor classes $a$ and $d$ on
$X$, the second and third Chern classes of $\F$ can be
written in the form
\[
c_2(\F) = a\,h+6\,h^2, \quad c_3(\F) = d\,h^2
\]
modulo higher codimension classes on $X$. The Grothendieck
Riemann-Roch theorem for the projection $p\: X\x P \to X$
easily gives (for example using \cite{schubert}) that
\[
-c_1(p_!\F(h)) = \frac52\, a- \frac12\, d.
\]
The locus $\xi\i B\sub X$ is set-theoretically the support
of $R^1 p_*\F(h)$. It is not hard to see that one can take
the family $X$ in such a way that the 0-th Fitting ideal
of $R^1 p_*\F(h)$ is actually reduced. Therefore the left hand
side of the equation above is $\xi^*\beta$.  On the other hand,
$a=\xi^*\alpha$ by the usual definition of the
$\mu$ map of Donaldson
\cite{Dona-1}, and $d=\xi^*\delta$.  Since the family
$\F/X$ was arbitrary, the required relation is actually
universal, and so holds also in $\Pic(M_6)\*\QQ$. (It suffices
to take a family with the properties that (i) $\xi^*\:\Pic_\QQ(M_6)
\to \Pic_\QQ(X)$ is injective, (ii) the Fitting ideal above is
reduced, and (iii) the general non-locally free sheaf in the
family has colength 1 in its double dual.)
\end{pf}

With this, we complete the proof of the second part of
\thmref{thm2} in the following way. The general fiber of
$f_6$ restricted to $\Delta$ has dimension 1, so
$f_6(\Delta)$ has dimension 19, see e.g.~\cite{Stro-1}.
Therefore we get
\[\begin{aligned}
\int_{H_7} s_{14}(\E\*\L) &= \int_{D_6}(\lambda+\tau)^{20} \\
&= \int_{M_6} \beta\, \alpha^{20} \\
&= \int_{M_6} (\frac52\, \alpha - \frac12\,\delta)\,\alpha^{20} \\
&= \frac52\int_{M_6} \alpha^{21} -\frac12\int_{\Delta}\alpha^{20}
= \frac52\, q_{21}.
\end{aligned}
\]

\section{A geometric interpretation}
\begin{defn}
A \emph{Darboux configuration} in $P\v$ consists of a
pair $(\Pi,C)$ where $\Pi\sub P\v$ is the union of $n+1$
distinct lines, no three concurrent, and $C\sub P\v$ is a
curve of degree $n$ passing through all the nodes of $\Pi$.
\end{defn}
If we let $Z\sub P$ consist of the $n+1$ points dual to
the components of $\Pi$, we have by Hulsbergen's theorem
\cite[thm.~4]{Bart-2} a natural 1-1
correspondence between Hulsbergen sequences
\eqref{Hulsbergen} and Darboux configurations $(\Pi,C)$,
by letting $C=J(F)$. Therefore $D_n$ can be used as a
compactification of the set of Darboux configurations, and
the intersection number
\[
\int_{D_n} \lambda^i (\tau+\lambda)^{3n+2-i} =
\int_{H_{n+1}} c_1(\L)^i s_{2n+2-i}(\E\*\L)
\]
can be interpreted as the number of Darboux configurations
$(\Pi,C)$ where $\Pi$ passes through $i$ given points and
$C$ passes through $3n+2-i$ given points.

It is not known whether the Barth map has degree 1 for
$n\ge5$. A related question is the following: Let
$(\Pi,C)$ be a general Darboux configuration ($n\ge 5$).
Is the inscribed polygon $\Pi$ uniquely determined by $C$?

\end{document}